\documentclass[12pt]{amsart}
\usepackage{
  graphicx,dbnsymb,amssymb,epic,eepic,picins,color,hyphenat,
  multicol,floatflt,slashbox
}
\usepackage{amsmath}
\usepackage[all]{xy}

\theoremstyle{plain}

\newtheorem{theorem}{Theorem}
\newtheorem{proposition}{Proposition}[section]

\newtheorem{lemma}[proposition]{Lemma}
\newtheorem{corollary}[proposition]{Corollary}

\theoremstyle{definition}
\newtheorem{definition}[proposition]{Definition}

\theoremstyle{remark}
\newtheorem{example}[proposition]{Example}

\newlength{\standardunitlength}
\catcode`\@=11
\long\def\@makecaption#1#2{%
    \vskip 10pt
    \setbox\@tempboxa\hbox{
      \small\sf{\bfcaptionfont #1. }\ignorespaces #2}%
    \ifdim \wd\@tempboxa >\captionwidth {%
        \rightskip=\@captionmargin\leftskip=\@captionmargin
        \unhbox\@tempboxa\par}%
      \else
        \hbox to\hsize{\hfil\box\@tempboxa\hfil}%
    \fi}
\font\bfcaptionfont=cmssbx10 scaled \magstephalf
\newdimen\@captionmargin\@captionmargin=2\parindent
\newdimen\captionwidth\captionwidth=\hsize
\catcode`\@=12

\def\qed{{\hfill\text{$\Box$}}}

\newlength{\globalparindent}
\def\proof{{\par\noindent {\em Proof. } }}

\def\bbZ{{\mathbb Z}}

\def\calA{{\mathcal A}}

\def\calJ{{\mathcal J}}
\def\calM{{\mathcal M}}
\def\calO{{\mathcal O}}
\def\calP{{\mathcal P}}
\def\calS{{\mathcal S}}
\def\calT{{\mathcal T}}

\def\calD{{\mathcal D}}

\def\calMo{{\mathcal M}^{(o)}}
\def\calMko{{\mathcal M}_k^{(o)}}

\newcommand{\Cobol}{{\mathcal Cob}^3_{o/l}}

\newcommand{\Kob}{\operatorname{Kob}}

\newcommand{\Kobh}{{\operatorname{Kob}_{/h}}}

\newcommand{\Mat}{\operatorname{Mat}}

\newcommand{\Mor}{\operatorname{Mor}}
\newcommand{\Obj}{\operatorname{Obj}}

\renewcommand{\qed}{~\hfill$\square$}

\begin{document}
\newdimen\captionwidth\captionwidth=\hsize

\title{The Jones polynomial and the planar algebra of alternating links}

\author{Hernando Burgos Soto
}
\address{ Department of Mathematics\\University of Toronto\\Toronto Ontario M5S 2E4\\Canada } \email{hburgoss@math.toronto.edu}

\date{
  First edition: Jul.{} 16, 2008.
  This edition: Feb.~19,~2009
}

\subjclass{57M25}
 \keywords{
  Alternating planar algebra,
  Coherently alternating element,
  Gravity information,
  Planar algebra,
  Rotation number,
  Skein Module,
  Smoothing,
  Tangles.
}

\thanks{This work was partially supported by Universidad del Norte in Barranquilla, Colombia, and Grant PRE0045002032, Colciencias-Unal 185, Colombia.\\
  \indent Electronic version: {\tt http//arxiv.org/PS\_cache/arxiv/pdf/0807/0807.2600.v1.pdf
  }
  .
}

\begin{abstract}
  It is a well known result from Thistlethwaite \cite{Thistlethwaite:SpanningTreeExpansion} that the Jones polynomial of a non-split alternating link is alternating. We find the right generalization of this result to the case of non-split alternating tangles. More specifically: the Jones polynomial of tangles is valued in a certain skein module, we describe an alternating condition on elements of this skein module, show that it is satisfied by the Jones invariant of the single crossing tangles $(\overcrossing)$ and
$(\undercrossing)$, and prove that it is preserved by appropriately ``alternating" planar algebra compositions. Hence, this condition is satisfied by the Jones polynomial of all alternating tangles. Finally, in the case of $0$-tangles, that is links, our condition is equivalent to simple alternation of the coefficients of the Jones polynomial.
\end{abstract}


\maketitle

\tableofcontents

\section{Introduction} \label{sec:intro}
Two years after the Jones polynomial was defined \cite{Jones:PolynomialInvariant}. M. B Thistlethwaite \cite{Thistlethwaite:SpanningTreeExpansion} based on W. T. Tutte's concepts of internal and external activities of edges with respect to a spanning tree and proved that the Jones polynomial of a non-split alternating link is an alternating polynomial. Subsequently, V. C. Turaev \cite{Turaev:JonesForTangles} presented an extension to tangles of the Jones polynomial which was a simple generalization of the Kauffman bracket presented in \cite{Kauffman:StatesModels}. Tangles are link pieces and the Jones invariant for tangles composes well, so it makes sense to ask what properties non-split alternating tangles have to produce (when they are composed) links with alternating Jones polynomials.\\
\indent In this paper we study the Jones polynomial of non-split alternating tangles and obtain a generalization of the Thistlethwaite result. For doing that, we define a special free $\mathbb{Z}[q,q^{-1}]$-module $\calMko$ with basis the set of oriented smoothings with $k$ strands and no circles. The orientation of the strands in an oriented smoothing $\sigma$ allows us to define a certain parameter associated with it, its {\it rotation number} $R(\sigma)$. We further define a subset of $\calMko$ formed by elements $P=\sum_{i=0}^nA_i\sigma_i$ where each $A_i$ is an alternating polynomial whose parity is related linearly with the rotation number of $\sigma_i$. The elements of this set whose partial closures have the same property, form what we call the collection of {\it coherently alternating elements}.This last collection will be denoted by $\calA$.\\
We introduce the concept of {\it alternating planar algebra} (an alternating planar algebra is a planar algebra that carries some extra restrictions). With these definitions we can state our two main theorems:
\begin{theorem}\label{theo:MainTheo1} The collection $\calA$  form an alternating planar algebra (that is, it is closed under compositions in alternating planar
diagrams).
\end{theorem}
 Our second theorem follows from the first; for it
reduces that proof to the simple task of verifying that the Jones polynomial of the one-crossing tangles $(\overcrossing)$ and
$(\undercrossing)$ (which are of course alternating) are
coherently alternating:
\begin{theorem}\label{theo:MainTheo2} Let $T$ be a non-split alternating $2k$-boundary tangle ($k>0$), then the Jones polynomial $\hat{J}(T)$ can be interpreted as an element of $\calA$.
\end{theorem}
It is a simple matter to verify that in the case of alternating
tangles with no boundary, i.e., in the case of alternating links,
this statement reduces to the Thistlethwaite result about alternating tangles.\\
\indent The work is organized as follows. In section
\ref{sec:Planar}, we review the Turaev construction for Jones polynomial of a tangle and the concept of planar algebra and present a few examples of planar algebras. Proposition \ref{prop:MorphismBracket}, introduced there will be a useful tool for the proof of theorem \ref{theo:MainTheo2}.
Section \ref{sec:Alternating} is devoted to introduce an
orientation in the smoothings, in this section we present the concept of rotation number and
some features of the $d$-input planar diagrams we are working with, which are used to define the concept of alternating planar algebras. We prove here that for an alternating planar algebra, the collection of operators is generated by operators falling into two classes, the {\it basic operators}. The proof of Theorem \ref{theo:MainTheo1} is based on this statement.
Section \ref{sec:On-Diagonal}  introduces the concepts of alternating elements in the skein module, coherently alternating elements and their partial closures. We state here some results about the elements obtained when a basic operator is applied to alternating elements, leading to the prove in section 4.5 of Theorem \ref{theo:MainTheo1}.. Finally in section
  \ref{sec:Examples} is dedicated to the study of non-split alternating tangles. We include special information in the strands of the tangles ({\it the gravity information}), which will give us not only an alternating orientation in the smoothings, but also a way to compose the tangles using alternating planar diagrams. Here, we prove Theorem \ref{theo:MainTheo2} and derive from it Thistlethwaite theorems formulated in \cite{Thistlethwaite:SpanningTreeExpansion}.\\
  \indent Analogously to the alternating property of the Jones polynomial, the reduced Khovanov homology of an alternating link is ``on diagonal". There ought to be a generalization of this property for tangles. In the spirit of this paper. We hope to address this question in a future publication.

\subsection{Acknowledgement} I wish to thank D. Bar-Natan, for many
helpful conversations we had at the University of Toronto. I would also
like to thank P. Lee, L. Leung and N. Martin for their comments and
suggestions.
\section{Jones polynomial and planar algebras, a review} \label{sec:Planar}
\subsection{The Jones polynomial for tangles}
\indent This construction is a simple generalization of the Kauffman model from the Jones polynomial of a link. Let $k$ be a non-negative integer and $S_k$ denote the set of all isotopy classes of tangle diagrams with $2k$ ends, having no crossing points and no closed components. An element of $S_k$ is called a {\it simple smoothing}. Turaev's construction \cite{Turaev:JonesForTangles} associates to each $2k$-end oriented tangle $T$ an element of the free $\mathrm{R}$-module $\mathrm{R}[S_k]$, where $\mathrm{R}=\mathbb{Z}[q,q^{-1}]$.\\
\indent A $k$-tangle diagram consists of $k$ open arcs and a finite number of loops in the disc $D^2$. The $2k$ end points of the arcs being $2k$ different fixed points in the boundary of $D^2$, and with the over/under crossing information represented as usual. We denote $M_k$ the $\mathrm{R}$-module generated by all ambient isotopic $k$-tangle diagrams modulo the relations:\\
 \begin{enumerate}
 \item The relation $D\cup \bigcirc=(q+q^{-1})D$, where $D$ is an arbitrary diagram and $\bigcirc$ is a loop.
 \item The so call skein relations \[ \overcrossing = q\smoothing - q^2\hsmoothing \hspace{2cm} \undercrossing = -q^{-2}\hsmoothing + q^{-1}\smoothing.
 \]
 \end{enumerate}
Following Turaev \cite{Turaev:Invariants3Manifolds} we call the module $M_k$ the skein module corresponding to $q\in \mathrm{R}$.
\indent Given a diagram $D$ of a tangle we can expand $D$ as a linear combination of elements in $S_k$ by applying relations $(1)$ and $(2)$, so $M_k$ is a free $\mathrm{R}$-module with a basis represented by simple smoothings. The element of $M_k$ representing the tangle diagram $D$ is an invariant of $T$; we call it the Jones polynomial of $T$ and denote it  $\hat{J}(T)$.
\subsection{Planar algebras}
The concept was introduced by Jones in \cite{Jones:PlanarAlgebrasI} to study the structure of subfactors, also see \cite{Jones:QuadraticTangles}.
  \parpic[r]{\includegraphics[scale=.45]%
{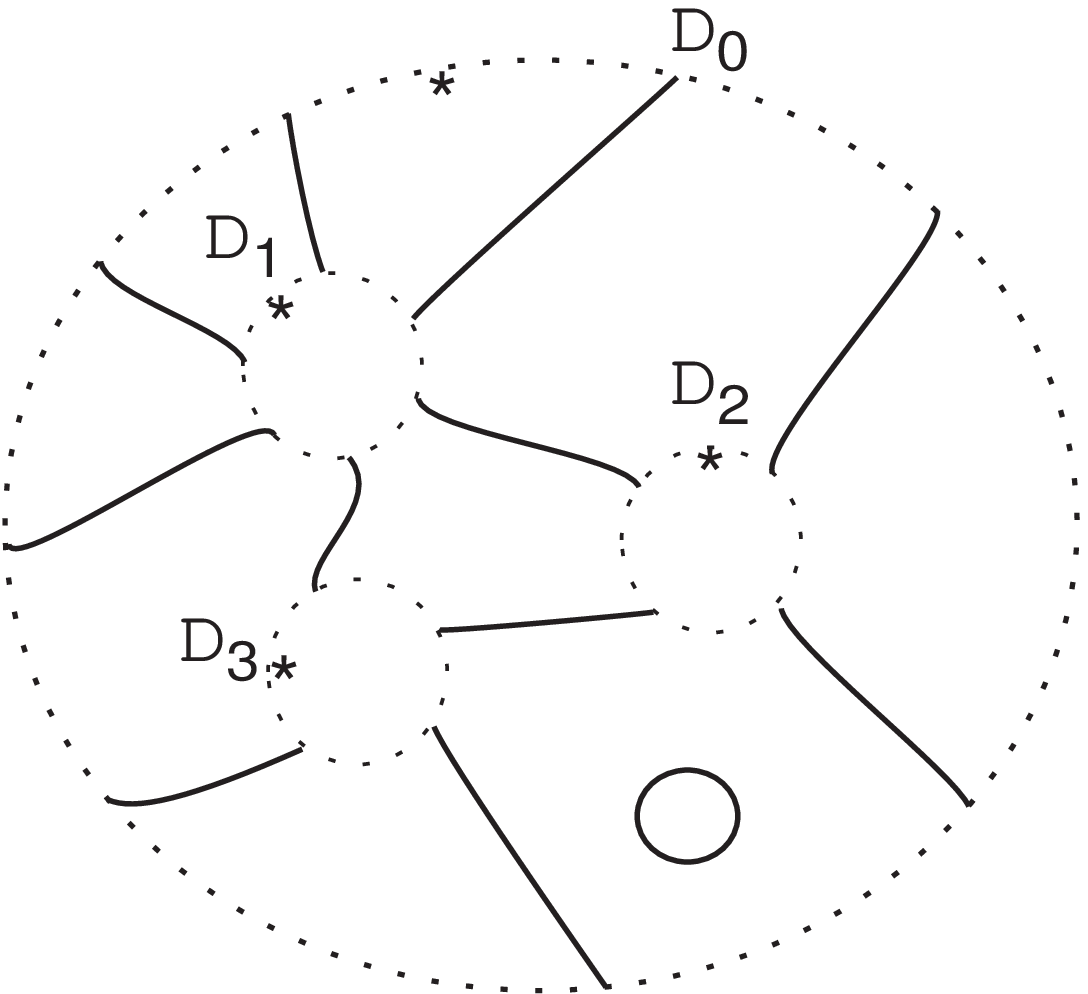}}Let $d$, $k_0,...,k_d$ be non-negative integer numbers. A \emph{d-input planar arc diagram} $D$ consists of a closed unit disc $D_0$, and discs $D_1,...,D_d$ smoothly embedded in its interior, a finite collection of circles that we will call loops, $\sum_{i=0}^dk_i$ disjoint open arcs ending transversally in the boundary of the discs, and a marked point for each of the disc. $2k_i$ points in the boundary of each disc $D_i$ are boundary points of the open arcs.  Arcs and loops are embedded in $int(D_0)\backslash\bigcup_{i=1}^dD_i$. The information is considered up to planar isotopy.
We say that $D$ is an oriented or an unoriented d-input planar arc diagram depending on whether its arcs are oriented or not.

\begin{figure}[hbt]\centering
\begin{tabular}{ccccc}
 \includegraphics[scale=.3]%
{figs/3planaralg.eps} & \hspace{.5cm} & \includegraphics[scale=.3]%
{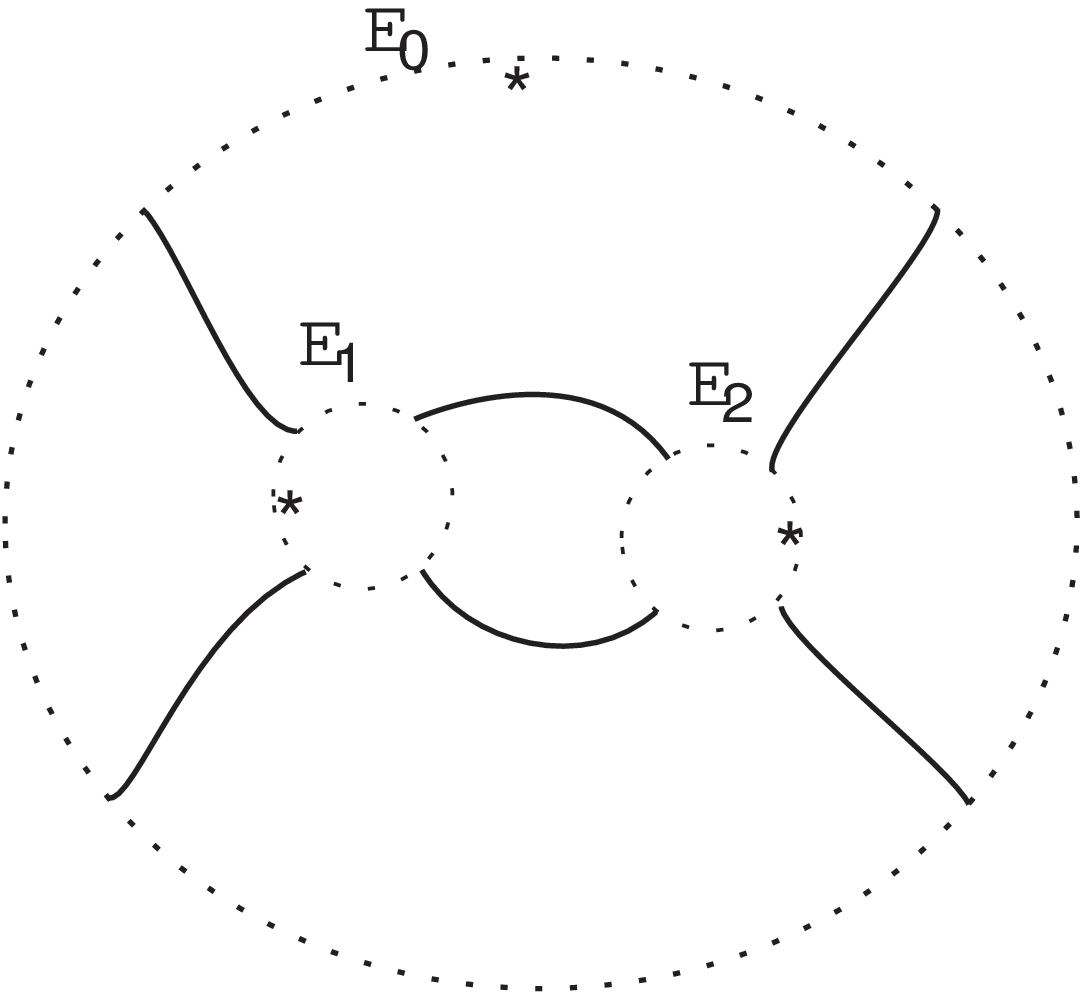} & \hspace{.5cm} &\includegraphics[scale=.3]%
{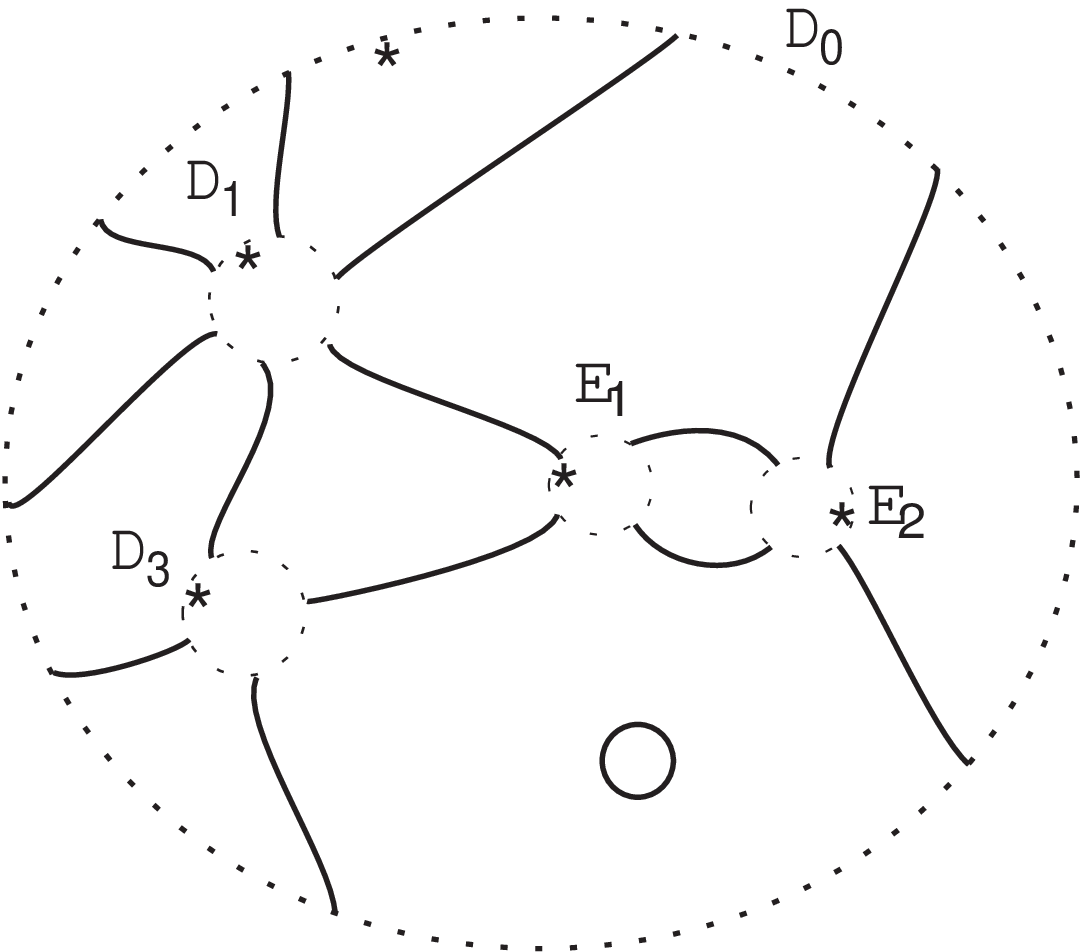}\\
$D$ & \hspace{.5cm} & $E$ & \hspace{.5cm} & $D\circ_2 E$
\end{tabular}
\caption{The composition of planar arc diagrams} \label{Fig:compoplanar}
\end{figure}

There is a natural way of composing planar arc diagrams. If $E$ is an $e$-input planar arc diagram, with $k$ arc boundary points in its boundary disc $E_0$, and $D$ is a $d$-input planar diagram with $k$ arc boundary points in its input disc $D_i$, then $D_i$ can be replaced in $D$ by $E$ by matching the corresponding marked points, removing the external disc of $E$ and connecting the arcs in a smoothly, producing in this manner a new $(d+e-1)$-input planar arc diagram. See Figure \ref{Fig:compoplanar}. We say that $E$ has been inserted in the i-th disc of $D$ and denote the new disc by $D\circ_i E$ .\\
  \indent In this way, a $k$-smoothing will be a $0$-input planar arc diagram with $k$ arcs ending on its external boundary. So, A $d$-input planar arc diagram defines a smoothing operation, i.e., by inserting appropriately $d$ smoothings in a $d$-input planar diagram we obtain a new smoothing. Actually, this is not the unique kind of object in which a planar arc diagram can define an operation, so we generalize the notion by the following\\
\indent Let $R$ be a ring. A planar algebra is a triple $\{\calP,\calD,\calO\}$ where $\calP$ is a sequence of $R$-modules: $\calP_0,\calP_1,...$, $\calD$ is a collection of input planar diagrams, and $\calO$ is a collection of operators. These operators are multilinear maps:  \[ O:\calP_{k_1}\times \cdots \times \calP_{k_d} \rightarrow \calP_{k} \] where $k,k_1,k_2,...$ are non-negative integers. There exists a bijection  between  $\calO$ and $\calD$, and a correspondence between the composition of planar arc diagrams and composition of operators. That is to say, if $O_D$ denotes the operator assigned to a $d$-input planar diagram $D$ then $O_{D\circ_i E}=O_D\circ_i O_E$. By abuse of notation and from now on, we shall denote the operator defined from a $d$-input planar diagram $D$ by the same symbol $D$.\\
\indent Given two planar algebras $\mathcal{P}^a$ and $\mathcal{P}^b$, a morphism of planar algebra is a collection of maps $\Phi_k:\mathcal{P}_k^a\rightarrow \mathcal{P}_k^a$ satisfying that for every d-input planar arc diagram $D$ with $2k_i$ arcs ending in the i-th input disc, we have that \begin{equation}\Phi_{k_0} \circ D= D\circ(\Phi_{k_1}\times \cdots\times\Phi_{k_d}) \end{equation}
    Here are some examples of planar algebras:
    \begin{itemize}
      \item The entire collection of tangles $\calT$: given a $d$-input planar diagram $D_0$ with input discs $D_1,...,D_d$, each $D_i$ with $2k_i$ marked points, then this planar arc diagram defines an operation $D:T_{k_1}\times \cdots\times T_{k_d}\rightarrow
          T_{k}$. In the same way, we can choose an oriented d-input arc planar diagram to define a similar operation using collections of oriented tangles.
      \item The entire collections associated to the Khovanov homology, $\Obj(\Cobol)$, $\Mor(\Cobol)$, $\Obj(\Mat(\Cobol))$, $\Mor(\Mat(\Cobol))$, $\Kob$ and $\Kobh$ are examples of planar algebras (see \cite[Section 5]{Bar-Natan:Cobordisms} for details)
      \item The collection of skein modules $\calM_k$ gives an example of an oriented planar algebra. Given a d-input oriented planar arc diagram with input discs $D_1,...,D_d$, each $D_i$ with $2k_i$ arc boundary points, the insertion of an element  $A_1\sigma_1+\cdots + A_n\sigma_n$ of $\calM_{k_i}$ in $D_i$ produces a linear combination of $(d-1)$-planar diagrams \[A_1D^{1i}+\cdots + A_nD^{ni}\] where $D^{ij}$ is the $(d-1)$-planar diagram resulting from putting $A_j$ in the input disc $D_i$ of $D_0$. It is clear now that we can insert in the same way an element of $\calM_{k_l}$ in the corresponding disc of $D^{ij}$ obtaining in this way a  linear combination of $(d-2)$-input planar arc diagrams. We can continue in this way until each of the $d$ initial input discs are filled, obtaining a linear combination of no crossing tangles. The terms in this linear combination could have some loops but each time that a loop appears we can just remove it and multiply the corresponding term by the factor $(q^{-1}+q)$. As we have a finite number of input disc in any planar diagram, it is clear that this process defines an operator as above.
    \end{itemize}
    The following proposition is an immediate consequence of the definitions.
    \begin{proposition}\label{prop:MorphismBracket} The Jones polynomial defines an oriented planar algebra morphism $\hat{J} :\calT \rightarrow \calM$.
    \end{proposition}
\section{Alternating oriented skein modules} \label{sec:Alternating}
We introduce an orientation on the elements of $\calS_k$. This
information will also allow us to state a special composition
between elements of $\calS_k$ and consequently among elements of $\calM_k$.
\subsection{Orientation of the smoothings}
The orientation in the smoothings is done in such a way that the orientation in the strands is alternating in the boundary of the disc where they are embedded. That can be achieved by shading the connected
components of the complement of the smoothing ({\it regions}) in the disc  black and
white, so that regions with common boundary have different shadings.
 This checkerboard coloring defines an alternating
orientation of the strings of the smoothing. We just have to draw arrows on the
strings so that the arrows point in a counterclockwise direction
around the black regions. We can forget the coloring and consider that the regions are divided in: {\it negative}, those with boundary oriented clockwise; and {\it positive}, those with boundary oriented counterclockwise.\\
\indent
\parpic[r]{\includegraphics[scale=.45]%
{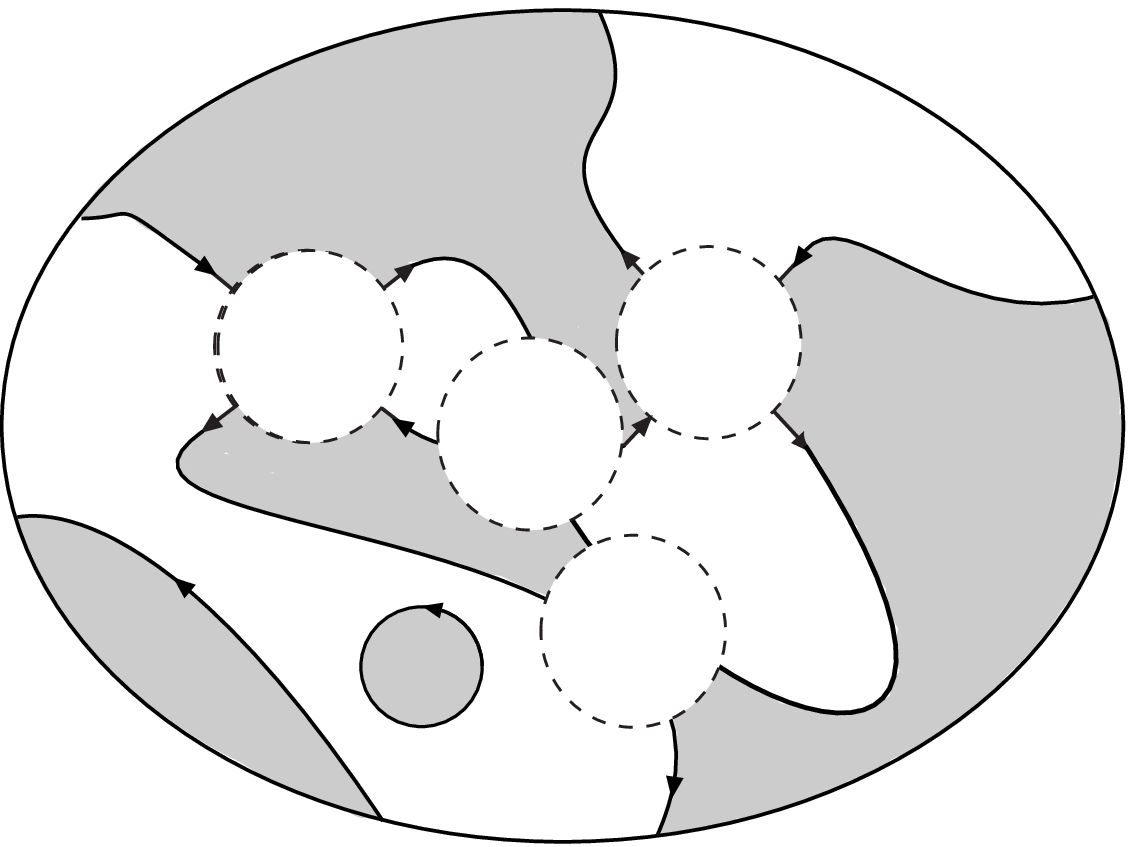}}After removing loops, the resulting collection of alternating oriented object obtained will be denoted with the symbols $\calS^{(o)}_k$ and the corresponding skein module $\calM_k^{(o)}$. A $d$-input
planar diagram with an alternating orientation of its arcs, which could be also achieved by a checkerboard coloring of the disc,  provides a good tool for the
horizontal composition of objects in $\calS^{(o)}_k$. Given smoothings
$\sigma_1,...,\sigma_d$, a suitable alternating $d$-input planar diagram $D$ to
compose them has the property that the $i$-th input disc has as many
boundary arc points as $\sigma_i$. Moreover placing $\sigma_i$ in the $i$-th input
disc, the orientation of $\sigma_i$ and $D$ match.

\subsection{Rotation numbers}
Given an oriented smoothing $\sigma$ possibly with loops, a strand in $\sigma$ is a 1-dimensional oriented manifold and a point in the boundary of $\sigma$ can be considered as an in-boundary point or an out-boundary point depending on the orientation of the strands in this point. We can enumerate the boundary points of the $\sigma$ from 0 to $2k-1$ starting from an in-boundary point of a strand, counting counterclockwise,
and finishing in the boundary point to the left of the mentioned
in-boundary point. A opened strand in an oriented smoothing $\sigma$ whose in-boundary and out-boundary points are numerated respectively $a$ and $b$, is denoted by $(a,b)$

\begin{figure}[hbt] \centering
\begin{tabular}{ccc}
\begin{tabular}{c} \includegraphics[scale=1.8]%
{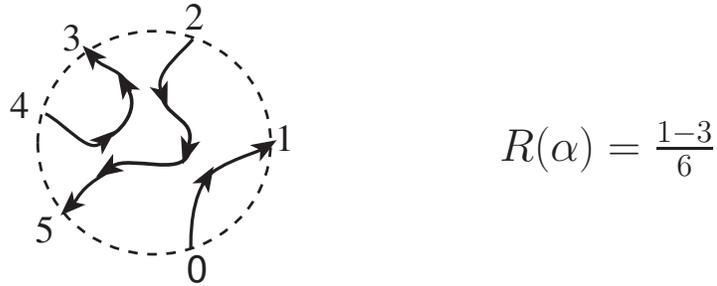} \end{tabular} & \hspace{1cm} & \begin{tabular}{c}
{\Large$R(\alpha)=\frac{1-3}{6}$}\end{tabular} \end{tabular}%
\caption{$\alpha =(0,1)$, $R(\alpha)=-\frac{2}{6}=-\frac{1}{3}$. The rotation number of
the complete resolution is 0}\label{Fig:rotate}
\end{figure}

\begin{definition} Given an open strand $\alpha$ of an alternating oriented smoothing $\sigma$, possibly with loops, enumerate the boundary points of $\sigma$ in such a way that $\alpha$ can be denoted by $(0,i)$. The {\it rotation number} of $(0,i)$, $R(\alpha)$, is
$\frac{i-k}{2k}$. If $\alpha$ is a loop, $R(\alpha)=1$ if $\alpha$ is oriented
counterclockwise, and $R(\alpha)=-1$ if $\alpha$ is oriented clockwise. The
rotation number of $\sigma$ is the sum of the
rotation numbers of its strings. See figure \ref{Fig:rotate}
\end{definition}

\begin{proposition}\label{prop:RotatingNumber} Let $(i_1,j_1)$ and $(i_2,j_2)$ be open strings bounding the same region of a smoothing $\sigma$, and let $\tau$ the smoothing resulting from replacing these two strings in $\sigma$ by the strings $(i_1,j_2)$ and $(i_2,j_1)$
then: (i) $R(\sigma)+1=R(\tau)$ if the region that the initial strings bound is positive, or (ii) $R(\sigma)-1=R(\tau)$ if the region that they bound is positive.
\end{proposition}
\begin{proof}
It is enough to observe every
string of $\sigma$, but the two involved in the change contribute to
the rotation number of $\tau$ in the same way. We just need to
observe how the rotation number of these two strands changes. Figure
\ref{Fig:negtoposopen} shows when the initial strings bound a positive region
strands.
\begin{figure}[hbt] \centering

\includegraphics[scale=.8]%
{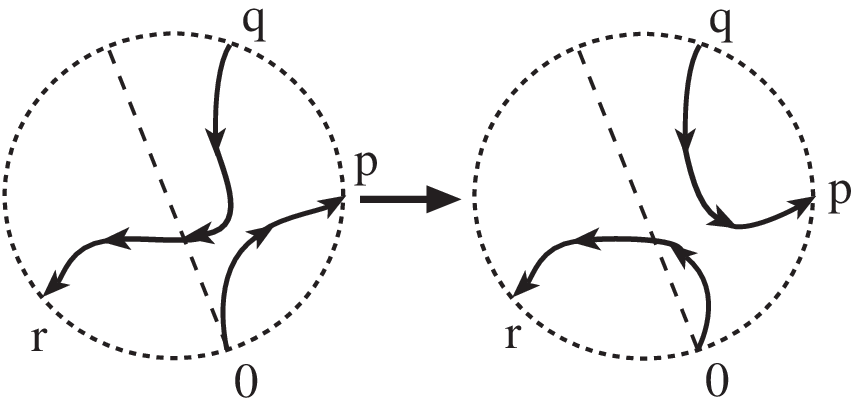}%
\caption{Change in the rotation of the smoothings by interchanging the end points of two strands}\label{Fig:negtoposopen}
\end{figure}
For a smoothing with $k$ strands, we can choose the entrance of one of
the strands as a 0 reference, if the boundaries of the strands  are
labeled with $0,p,q$ and $r$, and the initial arcs bound a positive region as in the figure. We can see that the
sum of the rotation numbers in the arcs of the initial smoothing is
$\frac{p-q+r-k}{2k}$ and
$\frac{p-q+r+k}{2k}$ in the resulting smoothing. \qed
\end{proof}
\begin{definition}A string in an alternating smoothing with $k$ ($k>0$) strands is called
minimal, if it is a negative loop or its rotation number is
$\frac{1-k}{2k}$, that is, it ends exactly to right of where it begins. It is
called maximal if it is a positive loop or its rotation number is
$\frac{k-1}{2k}$. A smoothing is minimal if all its strings are
minimal. It is called maximal if all its strings are
maximal.\end{definition}
\begin{figure}[hbt] \centering

\includegraphics[scale=.9]%
{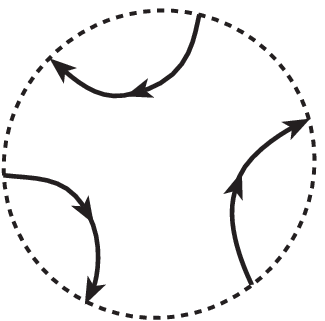} \hspace{1cm}
\includegraphics[scale=.9]%
{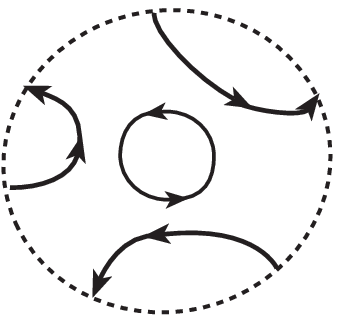}%
\caption{Examples of minimal and maximal
smoothings.}\label{Fig:maxmin}
\end{figure}
\indent Consider the checkerboard coloring of the disc. It is clear that every minimal smoothing is characterized by
having just one connected black region. If a minimal smoothing $\sigma_m$
with $k$ strands has $c_-$ negative loops, then the rotation number
of $\sigma_m$ is clearly $R(\sigma_m)=\frac{1-k}{2}-c_-$. In like manner, a
maximal smoothing is characterized by having only one white region.
If a maximal smoothing $\sigma_M$ with $k$ strands has $c_+$ positive
loops, then the rotation number of $\sigma_M$ is clearly
$R(\sigma_M)=\frac{k-1}{2}+c_+$.
\subsection{Connected alternating planar algebras}
\begin{definition} An alternating non-split planar diagram, see figure below,
 or type-$\mathcal{A}$ planar diagram $D$  is one with the
following properties:\begin{itemize}
\item the number $k$ of strings ending on the external boundary of
$D$ is greater than 0.
    \item There is complete connection among input discs of the
    diagram and its arcs, namely, the union of the diagram
    arcs and the boundary of the internal holes is a connected
    set.

    \item The arcs in the diagram are oriented and the in and
    out-strings alternate in every boundary component of the
    diagram.

\end{itemize}
\end{definition}
\parpic[r]{\includegraphics[scale=.45]%
{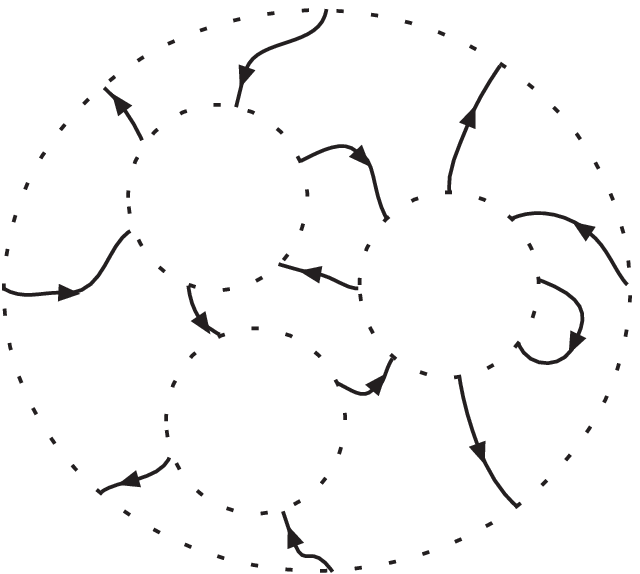}}A $d$-input type-$\mathcal{A}$ diagram  has an even number of
strings ending in each of its boundary components, and every string
that begins in the external boundary ends in a boundary of an
internal disk. Hence in these discs, we can classify the strings as:
{\it curls}, if they have its ends in the same input disc; {\it interconnecting arcs}, if its ends are in
different input discs, and {\it boundary arcs}, if they have one end in an input disc and the other in the external boundary
of the output disc. The arcs and the boundaries of the discs divide the surface of the diagram into disjoint regions. Some arcs and regions will be useful in the following definition and propositions.
\begin{definition}We assign the following numbers to every $d$-input planar diagram $D$:
\begin{itemize}
\item $i_D$: number of interconnecting arcs and curls, i.e., the number of non-boundary arcs.
\item $w_D$: number of negative internal regions. That is, in the checkerboard coloring, the white regions
whose boundary does not meet the external boundary of $D$.
\item $R_D$: the rotation associated number, which is given by the
formula \[ R_D=\frac{1}{2}(1+i_D-d)-w_D\]
\end{itemize}.
\end{definition}
Features of a $d$-input type-$\mathcal{A}$ planar diagram ensure
the following.

\begin{proposition}\label{prop:Diagramproperty}If $\sigma_1,\sigma_2,...,\sigma_d$ are $k_i$-strand smoothings, $k_i>0$, and $D$ is
type-$\mathcal{A}$ planar diagram, then $D(\sigma_1,\sigma_2,...,\sigma_d)$ is a
$k$-strand smoothing, $k=\sum_{i=1}^dk_i-i_D$. Furthermore, if
$\sigma_1,\sigma_2,...,\sigma_d$ are maximal (or minimal) smoothings and $D$
is a type-$\mathcal{A}$ planar diagram, then $D(\sigma_1,\sigma_2,...,\sigma_d)$ is
either a maximal (or a minimal) smoothing.
\end{proposition}
\begin{proof} Each of these statement follows immediately. For
example for the case in which every smoothing is minimal, all the
black regions in the big disc $D$ are connected for the connected
black region of the smoothings, so we obtain just one connected
black regions and hence a minimal smoothing.\qed
\end{proof}
\begin{proposition}\label{prop:AsoRotNumber}
Given the smoothings $\sigma_1,...,\sigma_d$ and a suitable d-input planar
diagram $D$, where every smoothing can be placed, the rotation
number of $D(\sigma_1,...,\sigma_d)$ is:
\begin{equation}\label{eq:AsocRotNumber}
R(D(\sigma_1,...,\sigma_d))=R_D+\sum_{i=1}^dR(\sigma_i)\end{equation}
\end{proposition}
\begin{proof}
Assume that $\sigma_1,...,\sigma_d$ are minimal smoothing, so by Proposition
\ref{prop:Diagramproperty} $D(\sigma_1,...,\sigma_d)$ is also a minimal
smoothing. Every internal white region will become a negative loop
in $D(\sigma_1,...,\sigma_d)$. If $2k$ is the number of non-internal strings in
$D(\sigma_1,...,\sigma_d)$, $R(D(\sigma_1,...,\sigma_d))$ will be given by
\[R(D(\sigma_1,...,\sigma_d))=\frac{1}{2}(1-k)+\sum_{i=1}^dc_{-i}-w_D\]
where $c_{-i}$ is the number of negative loops in $\sigma_i$.\\
\indent We have that $k=\sum_{i=1}^dk_i-i_D$ and
\[\sum_{i=1}^dR(\sigma_i)=\sum_{i=1}^d\left(\frac{1-k_i}{2}\right)-\sum_{i=1}^dc_{-i}\]
doing appropriate substitutions we obtain equation
(\ref{eq:AsocRotNumber}).\\
\indent For the general statement it is enough to observe that every
local change in a smoothing of one crossing in each $\sigma_i$ brings
the same local change in $D(\sigma_1,...,\sigma_d)$. That is, for example,
eliminating a loop in $\sigma_i$ brings the same change in $D(\sigma_1,...,\sigma_d)$ or interchanging end of the arcs bounded a positive region of $\sigma_i$ makes the same change or eliminates a negative loop in the whole picture, so the changes in
the rotation number in $D(\sigma_1,...,\sigma_d)$ is 1, which equals to the
change in the rotation number of $\sigma_i$. \qed
\end{proof}
\begin{definition} An alternating planar algebra is a triplet $\{\calP,\calD,\calO\}$ in which $\calP$, $\calD$, and $\calO$ have the same properties as in the definition of a planar algebra but with the collection $\calD$ containing only $\calA$-type planar diagrams.
\end{definition}
Diagrams with only one or two input discs deserves special attention. Operators defined from diagram like these are very important for our purposes since some of them are considered as the generators of the entire collection of operators in a connected alternating planar algebra.
\begin{figure}[hbt] \centering

\includegraphics[scale=.6]%
{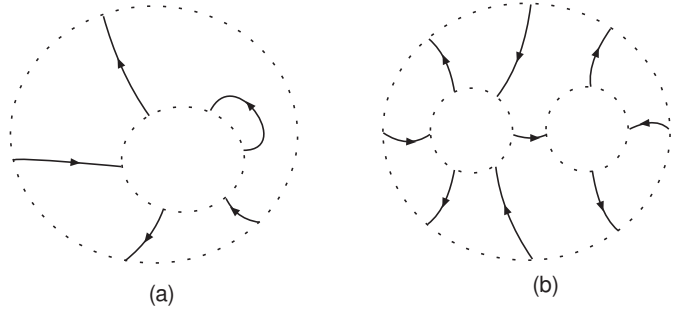}%
\caption{Examples of basic planar diagrams}\label{Fig:bplanar}
\end{figure}

\begin{definition} A basic planar diagram is a 1-input alternating planar
diagram with a curl in it,  or a 2-input alternating planar diagram with only one interconnecting arc. A basic operator is one defined from a basic planar diagram.
A negative unary basic operator is one defined from a basic 1-input diagram where the curl completes a negative loop. A positive unary basic operator is one defined from a basic 1-input diagram where the curl completes a positive loop. A binary operator is one defined from a basic 2-input planar diagram.
\end{definition}
\begin{proposition}\label{prop:RotNumberBasic}
The rotation associated number of a planar diagram belongs to
$\frac{1}{2}\mathbb{Z}$ and the case when we have a basic planar
diagram it is given as follows:
\begin{itemize}
\item If $D$ is  a negative unary basic operator, $R_D=-\frac{1}{2}$
\item If $D$ is  a  binary basic operator, $R_D=0$
\item If $D$ is  a positive unary basic operator, $R_D=\frac{1}{2}$
\end{itemize}
\end{proposition}
\proof{It follows from the definition of rotation associated number
and the features of the respective basic planar diagram. \qed}

\begin{proposition}\label{prop:compplanar} Any operator $D$ in an alternatively oriented planar algebra is the finite composition of basic operators.
\end{proposition}
\begin{proof} Every $d$-input $\calA$-type planar diagram is the composition of basic diagrams. \qed
\end{proof}

\section{The Planar Algebra $\calA$} \label{sec:On-Diagonal}
The coefficients in the terms of every element in $\calM_k^{(o)}$ are Laurent polynomial in $\bbZ$. Moreover we are interested in elements whose coefficients are alternating, so it will be useful to begin stating some concepts about alternating polynomials.
\subsection{Alternating polynomials}For our purpose, An alternating polynomial will be of the form $A(q)=\sum_{i=0}^{n}(-1)^{r+i}a_iq^{t+2i}$, where $n,r,t\in \bbZ$, and each $a_i$ is a non-negative integer. We say that two alternating polynomials $A(q)=\sum_{i=0}^{n}(-1)^{r+i}a_iq^{t+2i}$ and $B(q)=\sum_{i=0}^{n}(-1)^{r+i}b_iq^{t+2i}$ have the same parity if $a_ib_i>0$ \\
\indent Let $P=A_0\sigma_0+\cdots +A_m\sigma_m$ be an element of $\calM_k^{o}$, for $i\in \{0,...,m\}$, if $A_i$ is alternating it is possible to choose integers $r_i$, $t_i$ and $n_i$ in such a way that \begin{equation}\label{eq:PolyAlternating}A_i=(-1)^{r_i}q^{t_i}\sum_{j=0}^n(-1)^ja_{ij}q^{2j}\end{equation} with $a_{ij}\geq 0$, $a_{i0}\geq 0$ and $a_{in}\geq 0$. We call $(-1)^{r_i}a_{i0}q^{t_i}$ and $(-1)^{r_i+n}a_{in}q^{t_i+2n}$ respectively the minimal and maximal leading terms of $A_i$.
\subsection{Alternating elements}
\begin{definition} Given an element $P=A_0\sigma_0+\cdots +A_m\sigma_m$ of $\calMko$ with $\sigma_0$ a minimal smoothing, $\sigma_m$ a maximal smoothing, we say that $P$ is alternating if, for every $i\in \{1,...,m\}$, $A_i$ is an alternating polynomial and there exists $c_1,c_2\in \frac{1}{2}\bbZ$ such that the parity of $A_i$ is the same as the parity of $(-1)^{R(\sigma_i)+c_1}q^{R(\sigma_i)+c_2}$.
\end{definition}
Observe that the definition of alternating element of $\calMko$ ensures not only that the coefficients are alternating polynomials, but also that there is an alternation in the parity of these coefficients, i.e., when the rotation number increase in one, the exponents of $q$ in the coefficients change from even to odd or vice versa; when the rotation number increase by two, we obtain coefficients with the same exponents in $q$ but with opposite parity.
\begin{example}\label{Ex:alternatings} Here we have some examples of alternating elements in $\calMo$
\begin{enumerate}
\item $P_1=-q^{-2}\begin{tabular}{c}\includegraphics[scale=.3]%
{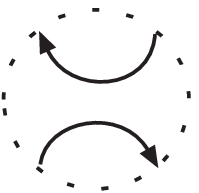}\end{tabular}+q^{-1}\begin{tabular}{c}\includegraphics[scale=.3]%
{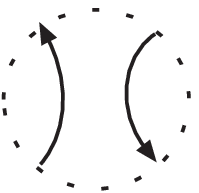}\end{tabular}$. This is the Jones invariant of the negative crossing $\undercrossing$, now with orientation in the smoothings. In this example the rotation number in the first term is $-\frac{1}{2}$ and in the second term it is $\frac{1}{2}$.
\item $P_2=(q^{-1}-2q)\begin{tabular}{c}\includegraphics[scale=.3]%
{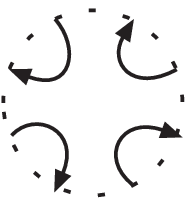}\end{tabular}+q^2\begin{tabular}{c}\includegraphics[scale=.3]%
{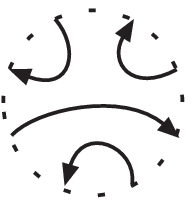}\end{tabular}+q^2\begin{tabular}{c}\includegraphics[scale=.3]%
{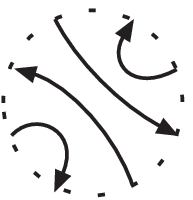}\end{tabular}+q^2\begin{tabular}{c}\includegraphics[scale=.3]%
{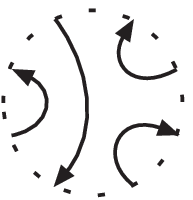}\end{tabular}-q^3\begin{tabular}{c}\includegraphics[scale=.3]%
{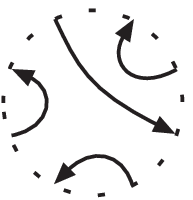}\end{tabular}$. Here the rotation number in the first term is $-\frac{3}{2}$; in the second, third and fourth term the rotation number is $-\frac{1}{2}$ and in the last term the rotation number is $\frac{1}{2}$.
\item $P_3=(q^{-2}-2+3q^2-2q^4)\begin{tabular}{c}\includegraphics[scale=.3]%
{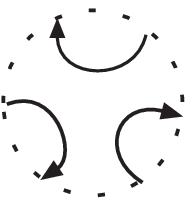}\end{tabular}+(-q^3+q^5)\begin{tabular}{c}\includegraphics[scale=.3]%
{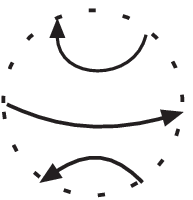}\end{tabular}+(q-q^3+q^5)\begin{tabular}{c}\includegraphics[scale=.3]%
{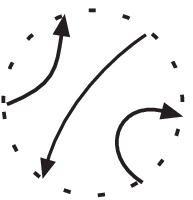}\end{tabular}+(-q^3+q^5)\begin{tabular}{c}\includegraphics[scale=.3]%
{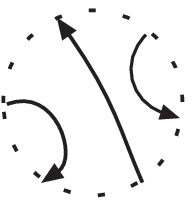}\end{tabular}+(q^4-q^6)\begin{tabular}{c}\includegraphics[scale=.3]%
{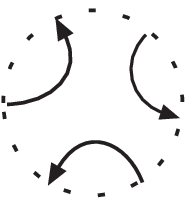}\end{tabular}$. In this last example, the rotation number in the first term is $-1$; in the second, third and fourth term the rotation number is $0$ and in the last term the rotation number is $1$.
\end{enumerate}
\end{example}

\subsection{Applying unary operators}
The elements of $\calMko$ can be inserted in appropriate unary basic planar diagrams obtaining in this way elements of $\calM_{k-1}^{(o)}$. This process can be summarized in the following steps:
\begin{enumerate}
\item placing of the element in the corresponding input disc of the unary basic planar diagram,
\item removing the loops obtained and replacing each of them by a factor $(q^{-1}+q)$, and
\item reducing the {\it like terms} resulting from the previous step.
\end{enumerate}
\begin{definition}Let $P=A_0\sigma_0+\cdots +A_n\sigma_n$ an element of $\calM_k$, then a partial closure of $P$ is an element $D_l\circ \cdots \circ D_1(P)$ where $0\leq l<k$ and every $D_i$ ($1\leq i\leq l$) is a unary basic operator.
\end{definition}
\parpic[r]{\includegraphics[scale=.4]%
{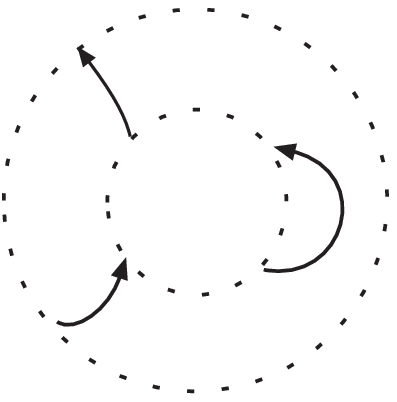}}The alternating element $P_1$ in the example \ref{Ex:alternatings} has the additional property that any of their partial closures are also alternating. For example embedding $P_1$ in a unary basic planar diagram $U_1$ as the one on the right, produce the element $U_1(P_1)=-q^{-2}\begin{tabular}{c}\includegraphics[scale=.3]%
{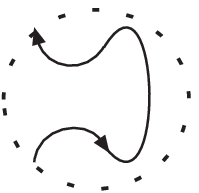}\end{tabular}+q^{-1}\begin{tabular}{c}\includegraphics[scale=.3]%
{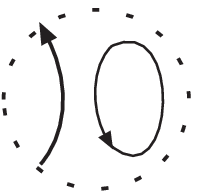}\end{tabular}=\begin{tabular}{c}\includegraphics[scale=.3]%
{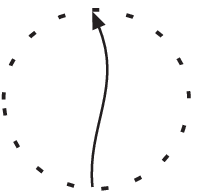}\end{tabular}$ which is obviously also alternating.
\begin{definition}
We say that an alternating element in $\calMko$ is {\it coherently alternating} if any of its partial closure is also alternating. We denote as $\calA_k$ the collection of all coherently alternating elements in $\calMko$. The symbol $\calA$ is used to denote $\bigcup_k\calA_k$.
\end{definition}
\parpic[r]{\includegraphics[scale=.4]%
{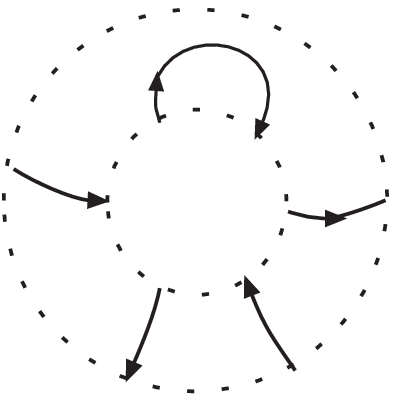}}Since the computation of other possible partial closures produces other alternating elements, the element $P_1$ of the example \ref{Ex:alternatings} is an element of $\calA_2$. Another example of coherently alternating element is the element $P_3$ of the same example, there are several possible partial closures, here we only calculate the one produced by inserting the element in the closure disc $C$ that appears on the right. It will be easy for the reader to compute the other partial closures. Inserting $P_3$ in $C$ produces the element \begin{eqnarray*}C(P_3)& = & (q^{-1}+q)(q^{-2}-2+3q^2-2q^4)\begin{tabular}{c}\includegraphics[scale=.3]%
{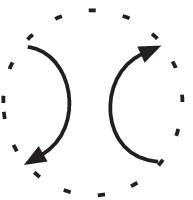}\end{tabular}+(q^{-1}+q)(-q^3+q^5)\begin{tabular}{c}\includegraphics[scale=.3]%
{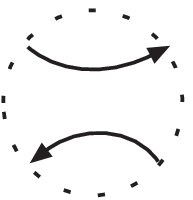}\end{tabular}\\
 &  & +(q-q^3+q^5)\begin{tabular}{c}\includegraphics[scale=.3]%
{figs/alternating11inv.eps}\end{tabular}+(-q^3+q^5)\begin{tabular}{c}\includegraphics[scale=.3]%
{figs/alternating11inv.eps}\end{tabular}+(q^4-q^6)\begin{tabular}{c}\includegraphics[scale=.3]%
{figs/alternating12inv.eps}\end{tabular}\\
 & = & (q^{-3}-q^{-1}+2q-q^3)\begin{tabular}{c}\includegraphics[scale=.3]%
{figs/alternating11inv.eps}\end{tabular}+(-q^2+q^4)\begin{tabular}{c}\includegraphics[scale=.3]%
{figs/alternating12inv.eps}\end{tabular}\end{eqnarray*}
which is also an alternating element.\\
\indent Let $P$ be an element of $\calMko$ and $C$ an appropriate closure operator such that $C(P)=\sum_{j=0}^nB_j\tau_j$. If $A_{i_1}\sigma_{i_1},...,A_{i_r}\sigma_{i_r}$ are terms in $P$ such that $C(A_{i_1}\sigma_{i_1}+...+A_{i_r}\sigma_{i_r})=B_r\tau_r$, we say that $A_{i_1}\sigma_{i_1},...,A_{i_r}\sigma_{i_r}$ are the generators of $B_r\tau_r$ under the closure $C$.\\
\indent Given a basic unary operator $U$ and an element $P$ of $\calMko$, for each term $A_i\sigma_i$ of $P$, $A_iU(\sigma_i)$ is not necessarily alternating. However we have the following proposition.
\begin{proposition}\label{prop:generatorterm} Let $P$ be an alternating term of $\calMko$ and let $U$ be an appropriate unary basic diagram such that $U(P)$ is alternating. Let $A_i\sigma_i$ and $A_j\sigma_j$ generators of the term $B\tau$ of $U(P)$, then:
\begin{enumerate}
\item If U is negative, then the parity of the minimal leading terms of the coefficient of $A_iU(\sigma_i)$  equals the parity of $B$.
\item If U is positive, then the parity of the maximal leading terms of the coefficient of $A_iU(\sigma_i)$ equals the parity of $B$.
\end{enumerate}
\end{proposition}

\begin{proof} Suppose $U$ is a negative basic planar diagram, inserting a smoothing in the diagram of $U$ can produce at most one loop, which is obviously negative. If there are only one generator for $B\tau$. Then this term must be alternating and both statements are proved. Assume, then, that $A_iU(\sigma_i)$ and $A_jU(\sigma_j)$ are two generators of $B\tau$. Clearly, after removing loops from $U(P)$, we have that the smoothing in $U(\sigma_i)$ and $U(\sigma_j)$ is $\tau$. However, before removing loops, there are exactly three possibilities for $U(\sigma_i)$ and $U(\sigma_j)$: (i) neither of them has the loop, (ii) both of them have the loop, and (iii) only one of them, say $U(\sigma_i)$, has the loop.\\
\indent If (i) holds, then by propositions \ref{prop:AsoRotNumber} and \ref{prop:RotNumberBasic} $R(\sigma_i)=R(\sigma_j)=R(\tau)+\frac{1}{2}$, and using the fact that $P$ is alternating, we conclude that $A_i$ and $A_j$ have the same parity, and of course the same happens with their minimal leading terms. Since neither $U(\sigma_i)$ nor $U(\sigma_j)$ has a loop, their minimal leading terms of the coefficients of $A_iU(\sigma_i)$ nor $A_jU(\sigma_j)$ are the same as $A_i$ and $A_j$.\\
\indent If (ii) holds, then by propositions \ref{prop:AsoRotNumber} and \ref{prop:RotNumberBasic} $R(\sigma_i)=R(\sigma_j)=R(\tau)-\frac{1}{2}$, so $A_i$ and $A_j$ have the same parity. Obviously, their minimal leading terms $(-1)^{r_i}a_{i0}q^{t_i}$ and $(-1)^{r_j}a_{j0}q^{t_j}$ have the same parity. After removing the common loop, the minimal leading term of the coefficients of $A_iU(\sigma_i)$ nor $A_jU(\sigma_j)$ are respectively $(-1)^{r_i}a_{i0}q^{t_i-1}$ and $(-1)^{r_j}a_{j0}q^{t_j-1}$ which also have the same parity.\\
\indent Finally, if $U(\sigma_i)$ has the loop and $U(\sigma_j)$ does not have it, then  $R(\sigma_i)=R(\sigma_j)-1$. Hence if the leading term of $A_i$ has the parity determined by $(-1)^{r}q^{t}$ then the minimal leading term of $A_j$ has the parity of $(-1)^{r+1}q^{t+1}$. Therefore, after removing the loops in $U(\sigma_i)$, the minimal leading term of the coefficients of $A_iU(\sigma_i)$ is $(-1)^{r}q^{t-1}$ which have the same parity as $(-1)^{r+1}q^{t+1}$.\\
\indent If $A_{i_1}\sigma_{i_1},...,A_{i_r}\sigma_{i_r}$ are the generators of $B\tau$, we have proved that the minimal leading term of the coefficients of $A_{i_1}U(\sigma_{i_1}),...,A_{i_r}U(\sigma_{i_r})$ have the same parity. Since one of these terms is the minimal leading term of $B$ the proof of the first statement is complete.\\
\indent The proof of the second statement is similar, we only need to change, in the former argument, minimal to maximal, $\frac{1}{2}$ to $-\frac{1}{2}$, and so on.\qed
\end{proof}
\subsection{Applying binary operators}
\begin{figure}[hbt] \centering

\includegraphics[scale=.5]%
{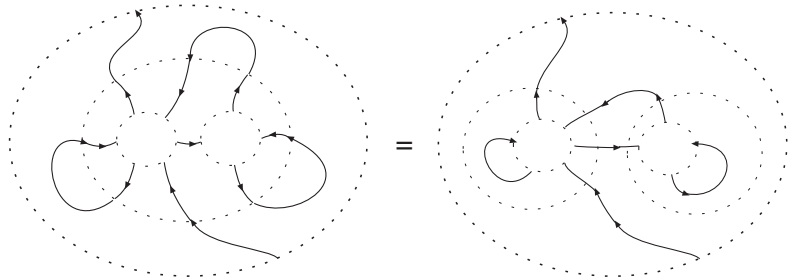}%
\caption{The closure of a binary operator can be considerate as the binary composition of two unary operator closures}\label{Fig:iqualplanars}
\end{figure}

\begin{proposition}\label{prop:ComBinaryOperator}
If $D$ is a suitable binary basic operator and $P,Q$ are alternating elements of $\calMko$, then $D(P,Q)$ is alternating.
\end{proposition}
\begin{proof}Let $P=\sum_{i=0}^{n}A_i\sigma_i$ and $Q=\sum_{j=0}^{m}B_j\tau_j$. Inserting $P$ and $Q$ in the disc $D$ produces an element $D(P,Q)=\sum_{i=0}^{n}\sum_{j=0}^{m}A_iB_jD(\sigma_i,\tau_j)$. Clearly, any $D(\sigma_i,\tau_j)$ is a smoothing with no loop, and by propositions \ref{prop:AsoRotNumber} and \ref{prop:RotNumberBasic} it has rotation number $R(D(\sigma_i,\tau_j))=R(\sigma_i)+R(\tau_j)$. If the minimal leading terms of $A_i$ and $B_j$ are respectively, $(-1)^{r_i}a_0q^{t_i}$ and $(-1)^{r_f}b_0q^{t_j}$, then $A_iB_j$ is an alternating polynomial with minimal leading term $(-1)^{r_i+r_j}a_0b_0q^{t_i+t_j}$. The relation between the exponents of $(-1)$ and $(q)$ in this last term with the rotation number of $D(\sigma_i,\tau_j)$ ensures that  $D(P,Q)$ is alternating.\qed
\end{proof}
\begin{proposition}\label{prop:clospolypoly} Let $P$ and $Q$ be elements of $\calA$ and let $D$ be a binary basic planar operator in which $D(P,Q)$ is well defined. For each partial closure $C(D(P,Q))$ there exists an operator $D'$ defined on a diagram without curls and elements $P',Q'$ in $\calA$ such that \[C(D(P,Q))=D'(P',Q')\]
\end{proposition}
\parpic[r]{\includegraphics[scale=.35]%
{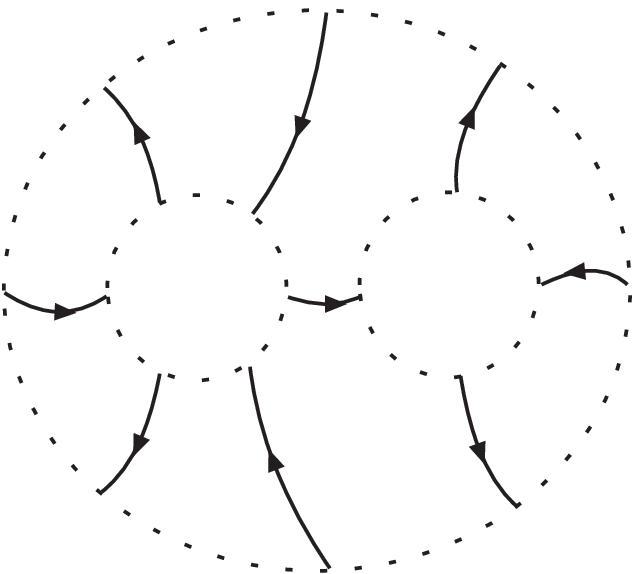}}
\begin{proof} \noindent We have a binary basic planar diagram $D$ as the one at the right. A closure of $D(P,Q)$ can be regarded as the composition of $P$ and $Q$ in an operator $C(D)$ defined from this closure, i.e., an operator defined from a disc $D$ embedded in a closure disc. Consider the strings with ends in the same input disc (the curls of the diagram). Since $D$ is a binary basic operator, in each input disc there is at least one string that is not a curl.  So we can regard the disc $D$ as a composition of two closure disc $E,E'$ embedded in a binary planar diagram $D'$ with no curl such that $C(D)=D(E,E')$. See Figure \ref{Fig:iqualplanars}.
 Hence $C(D(P,Q))=D(E(P),E'(Q))$. Since $E(P)$ and $E'(Q)$ are respectively closures of $P$ and $Q$ which are elements of $\calA$, the proposition is proved. \qed
\end{proof}

\begin{proposition}\label{prop:bimonomono}Let $\sigma$ and $\tau$ be elements of $\calS_k$ and $\calS_{k'}$ respectively, and let $D$ be a suitable binary planar operator defined from a no-curl planar arc diagram with output disc $D_0$, input discs $D_1,D_2$, and with at least one boundary arc ending in $D_1$, then there exists a closure operator $C$ and a unary operator $D'$ defined from a no-curl planar arc diagram such that $D(\sigma, \tau)=D'(C(\sigma))$. Moreover, if $P\in \calA_k$, then $D(P,\tau)$ is an alternating element of $\calMko$.
\end{proposition}
\begin{proof}Let $\calJ$ be the set of end points of the interconnecting arcs on the boundary of $D_2$, the disc in which $\tau$ is embedded. The strings in $\tau$ can be classified into three sets: $\Theta_i$, $i\in \{0,1,2\}$ formed by the arcs that have $i$ boundary points in $\calJ$. Since the points of $\calJ$ are consecutive points in the boundary of $D_2$, it is possible to embed a circle $C_0$ in $D(\sigma,\tau)$ with the following properties: \begin{enumerate}
\item$D_1$, and the elements of $\Theta_2$ are in the interior of $C_0$,
\item each arc of $\Theta_1$ is intersected once by $C_0$, and
 \item the elements of $\Theta_0$ are in the exterior of $C_0$. \end{enumerate} It is clear that the disc $C_0$ with the input disc $D_1$ defines now a closure operator $C$, and that $D(\sigma,\tau)\setminus C_0$ defines a no-curl unary operator $D'$ with input disc $C_0$. This prove the first statement.\\
\indent If $P=A_0\sigma_0+\cdots +A_n\sigma_n$, then  $D(P,\tau)=D'(A_0C(\sigma_0)+\cdots +A_nC(\sigma_n))=D'(C(P))$. The second statement follows from the fact that $P\in \calA$. \qed
\end{proof}

\subsection{ Proof of Theorem \ref{theo:MainTheo1}} Our first main result is an immediate consequence of the following statement.
\begin{lemma}\label{lem:colectofplanar} Let $P=A_0\sigma_0+\dots + A_n\sigma_n$ be a coherently alternating, and let $Q=B_0\tau_0+\cdots +B_m\tau_m$ be an alternating element of $\calMo$. Suppose that $D$ is a binary operator defined from a no-curl disc with at least one boundary arc ending in its first input disc,  then $D(P,Q)$ is an alternating element of $\calMo$.
\end{lemma}
\begin{proof}Assume that  $D$ has $l$ interconnecting arcs, we use induction over $l$. By proposition \ref{prop:ComBinaryOperator}, The lemma is true by $l=1$.\\
\indent Assume that if $D'$ has $l-1$ interconnecting arcs then $D'(P,Q)$ is alternating. A no-curl binary planar diagram with $l$ interconnecting arcs is obtained putting an $(l-1)$-interconnecting arc diagram $D'$, in a unary basic diagram $U$ in which the curl involves an arc of P. By proposition \ref{prop:bimonomono}, for each $j\in{1,...,m}$,  \[P'_j=\sum_{i=1}^nA_iB_jU\circ D'(\sigma_i,\tau_j)=\sum_{i=1}^{n_j}A'_{ij}\sigma'_{ij}\] is alternating.\\
\indent We have that before combining like terms $D(P,Q)=U\circ D'(P,Q)$ is the sum of alternating elements $P'_j$. To prove that after combining like terms in $D(P,Q)$ we obtain an alternating element, it is enough to show that given terms $A'_i\sigma'_i$ in $P'_i$ and  $A'_{ij}\sigma'_j$ in $P'_j$ such that $R(\sigma'_i)=R(\sigma'_j)$, $A'_i$ and $A'_j$ have the same parity.\\
\indent Let $\{A'_{i1}\sigma'_{i1},...,A'_{im_i}\sigma'_{im_i}\}$ and $\{A'_{i1}\sigma'_{i1},...,A'_{im_i}\sigma'_{im_i}\}$ be the respective set of generators of $A'_i\sigma'_i$ and $A'_j\sigma'_j$ under the closure $U$. By proposition \ref{prop:generatorterm} each minimal leading term in the coefficients of the elements of $G_i=\{A'_{ir_i}U(\sigma'_{ir_i})\}_{r_i\in \{1,...m_i\}}$ has the same parity as $A'_i$ and each minimal leading term in the coefficients of the elements of $G_j=\{A'_{jr_j}U(\sigma'_{jr_j})\}_{r_j\in \{1,...m_j\}}$ has the same parity as $A'_j$.\\
\indent Suppose that $U$ is a negative basic planar diagram and $A'_{ir_i}\sigma'_{ir_i}$ and $A'_{jr_j}\sigma'_{jr_j}$ are generators of $A'_i\sigma'_i$ and $A'_j\sigma'_j$ respectively. We just need to prove that the minimal leading terms in the elements of $G_i$ have the same parity as the minimal leading terms of the elements in $G_j$. The proof of this last statement is very similar to proposition \ref{prop:generatorterm}, i.e., before removing the loops in $D(P,Q)$, we have three options for $A'_{ir_i}U(\sigma'_{r_i})$ and $A'_{jr_j}U(\sigma'_{jr_j})$: (i) neither of them has the loop, (ii) both of them have the loop and (iii) just one of them have the loop. If neither of them has the loop , then by propositions \ref{prop:AsoRotNumber} and \ref{prop:RotNumberBasic}, $R(\sigma'_{ir_i})=R(\sigma'_i)+\frac{1}{2}=R(\sigma'_j)+\frac{1}{2}=R(\sigma'_{jr_j})$. Since $D'(P,Q)$ is alternating, we have that $A'_{ir_i}$ and $A'_{jr_j}$ have the same parity, and of course the same happens with their minimal leading terms. The remainder of the proof follows that of proposition \ref{prop:generatorterm}, taking into account that $D'(P,Q)$ is alternating. \qed
\end{proof}
\hspace{.5cm}
\begin{proof}(Of theorem \ref{theo:MainTheo1}.) By proposition \ref{prop:compplanar} we just need to prove that $\calA$ is closed under composition of basic operators. Let $P\in \calA$ and let $U$ be a basic unary operator. Since $U(P)$ is a partial closure of $P$, $U(P)$ is alternating. Furthermore any partial closure of $U(P)$ is also a partial closure of $P$, so $U(P)\in \calA$.\\
Let $P$ and $Q$ be elements of $\calA$ and $D$ a basic binary operator, Since $D$ is defined from a type-$\calA$ diagram, there is at least one boundary arc. Without loss of generality, we can assume that there is one boundary arc ending in the first input disc of $D$. By proposition \ref{prop:ComBinaryOperator}, $D(P,Q)$ is alternating. Let $C(D(P,Q))$ be a partial closure of $D(P,Q)$, by proposition \ref{prop:clospolypoly} there exist $P',Q'\in \calA$ and a binary operator $D'$ defined from a no-curl planar diagram such that $C(D(P,Q))=D'(P',Q')$. The fact that $D'(P',Q')$ is alternating follows immediately from lemma \ref{lem:colectofplanar}.  \qed
\end{proof}

\section{Non-split alternating tangles} \label{sec:Examples}
\subsection{Gravity information}
\begin{figure}[hbt]
\begin{tabular}{ccc}
\includegraphics[scale=0.6]%
{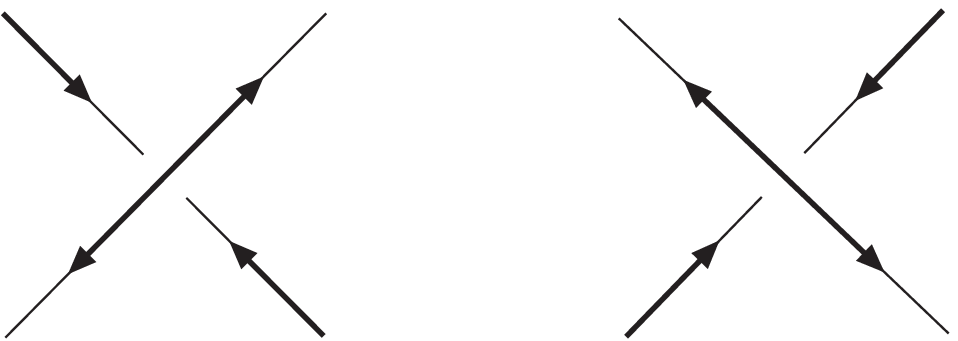} & \hspace{2cm} &\includegraphics[scale=.8]%
{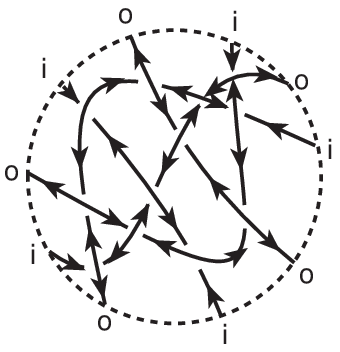} \\
(a) & & (b)
\end{tabular}
\caption{(a) Gravity information in a neighborhood of a crossing. (b) Gravity information in the tangle. We use: o for out-strands
and i for in-strands.
}\label{gravicross}
\end{figure}

Given a diagram of an alternating tangle we add to it some special
information which will help us to compose the Jones invariant
of an alternating tangle in an alternating planar diagram. That information is illustrated by drawing, in every
strand of the diagram, an arrow pointing in to the undercrossing, or
equivalently (if we have alternation), pointing out from the overcrossing. In a
neighborhood of a crossing the diagram looks like the one in Figure
\ref{gravicross}(a).
Figure \ref{gravicross}(b) shows a diagram of a tangle in which we have
added the gravity information to the whole tangle. Once it is done,
we find that the strand which meet the boundary of the disc, where
the tangle is embedded, are of two types. The ones that point
outside of the disc, which will be called {\it out-strands}, and
those which point inside of the disk, {\it in-strands}.\\
\indent We can observe, (see Figure \ref{gravicross}(a)) that if we make a
smoothing in the crossing, the orientation provided by the gravity
information is preserved, and that a $0$-smoothing is clockwise and
$1$-smoothing is counterclockwise, see figure \ref{Fig:smoothing}.
Also it is easily observed that if we go into a non-split
alternating for an in-strand and turn to the right (a
$0$-smoothing) every time that we meet a crossing, we are going to
get out of the tangle along the strand immediately to the right. So the in
and out-strands of the diagram of the tangle are arranged
alternatingly. These two observations are stated in the following
two propositions:
\begin{figure}[hbt] \centering

\includegraphics[scale=.8]%
{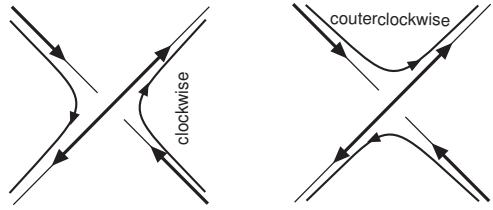}%
\caption{The smoothings in the diagrams preserve the gravity
information}\label{Fig:smoothing}
\end{figure}
\begin{proposition}\label{prop:gravit.preserv}The 0-smoothings and 1-smoothings preserve the gravity information, the first ones provides a clockwise orientation of the pair of strands in the smoothing, the
last provide a counterclockwise orientation.\end{proposition}
\begin{proposition}\label{prop:AlternanSmoothly}In any non-split alternating tangle, if the $k$-th strand is an in-strand the $(k+1)$-th strand is an out-strand.\end{proposition}
Propositions \ref{prop:gravit.preserv} and \ref{prop:AlternanSmoothly} leads to that
the smoothings of a tangle could be drawn as a trivial tangles in
which arcs are oriented alternatingly. So the Jones polynomial in fact produce an alternating planar algebra morphism.
\subsection{Proof of Theorem \ref{theo:MainTheo2}}
This proof is a direct application of proposition \ref{prop:MorphismBracket} applied to alternating planar algebras and Theorem \ref{theo:MainTheo1}.\\

\begin{proof} (Of Theorem \ref{theo:MainTheo2}) The Jones polynomial
of a 1-crossing tangle is a coherently alternating element. See the first example in
\ref{Ex:alternatings}. Any non-split alternating $k$-strand tangle $T$
with $n$ crossing, is obtained by a composition of $n$ of these
1-crossing tangles, $T_1$,...,$T_n$, in an $n$-input
type-$\mathcal{A}$ planar diagram. By proposition \ref{prop:MorphismBracket},
using the same $n$-input planar diagram for composing
$\hat{J}(T_1)$,...,$\hat{J}(T_n )$ we obtain the Jones polynomial of the original tangle. By Theorem \ref{theo:MainTheo1}, this is an element of $\calA$.  \qed \end{proof}
\subsection{Non-split alternating links}
It follows immediately from Theorem \ref{theo:MainTheo2} that
\begin{corollary}\label{col:1StrandTangle}The Jones polynomial $\hat{J}(T)$ of a non-split alternating tangle $T$ with only one open arc is a one term element in $\calMko$ whose smoothing is just an arc, and whose coefficient is an alternating polynomial in $\mathbb{Z}$.
\end{corollary}
\begin{corollary}The Jones polynomial $J(L)$ of a non-split alternating Link $L$ is an alternating polynomial in $\mathbb{Z}$.
\end{corollary}
\begin{proof}
Let $L$ be a non-split alternating link. $L$ is obtained by putting a
1-strand tangle $T$ in a 1-input planar arc diagram with one curl and no boundary arcs. So, we just have to put the Jones polynomial of this
1-strand tangle in the same 1-input planar arc diagram, so by corollary \ref{col:1StrandTangle} the ``unnormalized" Jones polynomial $\hat{J}(L)$ is a one term element of $\calMko$ whose smoothing is a loop and whose coefficient is an alternating polynomial with coefficients in $\mathbb{Z}$. To obtain the normalized Jones polynomial $J(L)$ of
the link $L$, we divide $\hat{J}(L)$ by $q+q^{-1}$ which is exactly the value of the circle obtained.  \qed
\end{proof}

\end{document}